\UseRawInputEncoding
\documentclass[12pt, reqno]{amsart}
\usepackage{amssymb}
\usepackage{amsfonts}
\usepackage{amssymb,amsmath,amscd}
\usepackage[colorlinks=true]{hyperref}
\usepackage{mathrsfs}
\newtheorem{Theorem}{\indent Theorem}[section]

\newtheorem{Lemma}[Theorem]{\indent Lemma}
\newtheorem{Corollary}[Theorem]{\indent Corollary}

\theoremstyle{remark}
\newtheorem{Remark}{Remark}

\begin{document}
\centerline{\bf On general sums involving the floor function with applications to $k$-free numbers }
 \pdfoutput=1
\bigskip
\centerline{\small Wei Zhang}
\bigskip

\textbf{Abstract}\
In this paper, we consider  sums related to the floor function.
We can improve previous results for some special arithmetic functions considered by   Bordell\`{e}s \cite{B}, Stucky \cite{St} and Liu-Wu-Yang \cite{Wu}.
It is worth emphasizing that we use much simpler methods to give much better results than previous.



\bigskip
\bigskip
\numberwithin{equation}{section}

\section{Introduction}
Recently, the sum
\[
S_{f}(x)=\sum_{n\leq x}f\left(\left[\frac{x}{n}\right]\right)
\]
has attracted  many experts' special attention (for example, see \cite{BDHPS,B,LWY,Wu,Zhai}),
where $f$ is a complex-valued arithmetic function and $[\cdot]$ denotes the floor function (i.e. the greatest integer function). One can call $S_{f}(x)$ the fractional sum of $f$ (see \cite{St}).

Specially, for some fixed $\eta\in(0,1)$ and
\[
f(n)\ll n^{\eta},
\]
independently, Wu \cite{Wu} and Zhai \cite{Zhai} showed that
\[
S_{f}(x)=C_{f}x+O\left(x^{(1+\eta)/2}\right),
\]
where $f$ is a complex-valued arithmetic function and
\[
C_{f}=\sum_{n=1}^{\infty}\frac{f(n)}{n(n+1)}.
\]
This formula improves the recent result obtained by
Bordell\`{e}s, Dai, Heyman, Pan and Shparlinski \cite{BDHPS}.

On the other hand, for some fixed $\eta\in(0,2)$ and
\[
\sum_{n\leq x}|f(n)|^{2}\ll x^{\eta},
\]
in \cite{BDHPS}, Bordell\`{e}s, Dai, Heyman, Pan and Shparlinski  proved that
\begin{align}\label{B}
S_{f}(x)=C_{f}x+O\left(x^{(1+\eta)/3}(\log x)^{(1+\eta)(2+\varepsilon_{2}(x))/6}\right),
\end{align}
where $f$ is a complex-valued arithmetic function,
\[
\varepsilon_{2}(x)=\left(\frac{2\log \log \log x}{\log \log x}\right)^{1/2}\left(1+\frac{30}{\log \log \log x}\right)
\]
 and
\[
C_{f}=\sum_{n=1}^{\infty}\frac{f(n)}{n(n+1)}.
\]
Recently,
for some fixed $\eta\in(0,2)$ and
\[
\sum_{n\leq x}|f(n)|^{2}\ll x^{\eta},
\]
 Wu  \cite{Wu} showed that
\begin{align}\label{Wu}
S_{f}(x)=C_{f}x+O\left(x^{(1+\eta)/3}(\log x)^{(1+\eta)\varepsilon_{2}(x)/6}\right),
\end{align}
where here and throughout, $C_{f}$ denotes the constant
\[
C_{f}=\sum_{n=1}^{\infty}\frac{f(n)}
{n(n+1)},
\]
$f$ is a complex-valued arithmetic function and
\[
\varepsilon_{2}(x)=\left(\frac{2\log \log \log x}{\log \log x}\right)^{1/2}\left(1+\frac{30}{\log \log \log x}\right).
\]
This formula improves the recent result obtained by
Bordell\`{e}s, Dai, Heyman, Pan and Shparlinski \cite{BDHPS} for the aspect of powers of $\log x$.

Recently, in \cite{ZW}, Zhao and Wu showed that
\begin{align}\label{zw0}
S_{f}(x)=
C_{f}x+O\left(x^{(2+3\eta)/8+\varepsilon}\right),
\end{align}
by using the following  result of \cite{Ba}
\begin{align}\label{ba0}
\sum_{n\geq1}\left(\left\{\frac{x}{n+1}\right\}
-\left\{\frac{x}{n}\right\}  \right)^{2}
=\frac{\zeta(3/2)\sqrt{x}}{\pi}+O(x^{3/7}).
\end{align}
We point out that  combining the main result of Balazard \cite{Ba}, we can obtain a slight better result.

The first aim of this paper is to show the following.
We can improve previous results and give non-trivial estimate under the assumption
\[
\sum_{n\leq x}|f(n)|^{2}\ll x^{\eta},
\]
with any fixed $\eta\in(0,2).$ We can improve   (\ref{zw0}) by eliminating the $x^{\varepsilon}$ term. Moreover, the proofs of us are much more elementary. The main idea is a basic observation such that
\[
\sum_{n\leq x}|f(n)|^{2}\ll x^{\eta}
\]
implies   $|f(n)| \ll x^{\eta/2}.$ Hence we can avoid using the so-called $r$-th Hooley divisor function used in \cite{BDHPS,Wu,ZW}.
\begin{Theorem}\label{th1}
Let $f$ be a complex-valued arithmetic function.
 Assume that
\[
\sum_{n\leq x}|f(n)|^{2}\ll x^{\eta},
\]
with any fixed $\eta\in(0,2).$ Then we have
\[
S_{f}(x)=
C_{f}x+O\left(x^{(2+3\eta)/8}\right),
\]
where
\[
C_{f}=\sum_{n=1}^{\infty}\frac{f(n)}{n(n+1)}.
\]
\end{Theorem}

In  general, for $f$  a positive real-valued arithmetic function,   one can obtain some much better results by using the theory of Fourier series and exponential sums. For example, one can refer to  \cite{B,LWY,MS,St}.
In this paper, we consider the following general result. We can improve previous results for some special arithmetic functions considered by   Bordell\`{e}s \cite{B}, Stucky \cite{St} and Liu-Wu-Yang \cite{Wu}.

\begin{Theorem}\label{th2}
Let $(\kappa,\lambda)\neq(0,1)$ be an exponent pair.
Let $f$ be a   positive real-valued arithmetic function such that  $f(n)=\sum_{d|n}g(d),$
and for any sufficiently large $x$
\[
\sum_{n\leq x}|g(n)|\ll x^{\lambda/(1+\kappa)}.
\]
Then we have
\begin{align*}
S_{f}(x)=
C_{f}x+
 \begin{cases}
O\left(x^{\lambda/(1+\lambda)}(\log x)^{\alpha}\right)\ \ &\textup{if}\ f(n)\ll 1,\\
\\
O\left(x^{\lambda/(1+\lambda)+\varepsilon}\right)\ \ &\textup{if}\ f(n)\ll n^{\varepsilon},
\end{cases}
\end{align*}
where
\[
C_{f}=\sum_{n=1}^{\infty}
\frac{f(n)}{n(n+1)}.
\]
and $\alpha=1$ if $(\kappa,\lambda)=(1/2,1/2)$ and 0 otherwise.
\end{Theorem}
Let $\mu_{k}(n)$ be the indicator function of the $k$-free numbers. Bordell\`{e}s \cite{B} considered the sum
\[
 S_{\mu_{2}}(x)=\sum_{n\leq x}\mu_{2}\left(\left[\frac{x}{n}\right]\right)
\]
and proved that
  \begin{align}\label{B1}
 S_{\mu_{2}}(x)=
\sum_{n=1}^{\infty}
\frac{\mu_{2}(n)}{n(n+1)}x
+O\left(x^{1919/4268+\varepsilon}\right).
 \end{align}
Recently, Liu-Wu-Yang \cite{LWY} improved the result of Bordell\`{e}s \cite{B} by showing that
\begin{align}\label{W1}
 S_{\mu_{2}}(x)=
\sum_{n=1}^{\infty}
\frac{\mu_{2}(n)}{n(n+1)}x
+O\left(x^{2/5+\varepsilon}\right).
\end{align}
In \cite{St}, for $k\geq3,$ Stucky remarked that
\begin{align}\label{S1}
 S_{\mu_{k}}(x)=\sum_{n\leq x}\mu_{k}\left(\left[\frac{x}{n}\right]\right)=
\sum_{n=1}^{\infty}
\frac{\mu_{k}(n)}{n(n+1)}x
+O\left(x^{\theta_{k}}\right),
\end{align}
 where
 \[
 \theta_{k}=\left(1+\frac1k\right)
 \left(3+\frac1k\right)^{-1}.
 \]
By using Theorem \ref{th2}, we can obtain the following corollary by choosing suitable exponent pairs.
\begin{Corollary}\label{co}
We have
\begin{align*}
 S_{\mu_{k}}(x)=
 \begin{cases}
C_{\mu_{k}}x+O\left(x^{11/29}(\log x)^{2}\right)\ \ &\textup{if}\ k=2,\\
\\
C_{\mu_{k}}x+O\left(x^{1/3}(\log x)\right)\ \ &\textup{if}\ k\geq3,\\
\end{cases}
\end{align*}
where
\[
C_{\mu_{k}}=\sum_{n=1}^{\infty}\frac{\mu_{k}(n)}{n(n+1)}
\]
and $\mu_{k}(n)$ is the indicator function of the $k$-free numbers.
\end{Corollary}
One can find that $11/29<0.3794$ and $2/5=0.4$. Hence according to (\ref{B1}) and (\ref{W1}), we can give much better results. For $k\geq3$,  we find that
\[
 \theta_{k}=\left(1+\frac1k\right)
 \left(3+\frac1k\right)^{-1}>1/3.
\]
Hence, our result is   better than (\ref{S1}) for fixed $k\geq2.$

On the other hand, Stucky \cite{St} also considered the fractional sum of $f(n)=\sum_{d|n}\frac{1}{d^{\beta}}.$
For $\beta\in(2/3,1]$ (the result is trivial for $0<\beta\leq 2/3$, see \cite{St} for details), it is proved that
\[
S_{f}(x)=C_{f}x+O(x^{\frac{2-\beta}{2+\beta}}).
\]
For this special situation, we can obtain the following result.
\begin{Corollary}
Let
$$f(n)=\sum_{d|n}\frac{1}{d^{\beta}}.$$
For $\beta\in(2/3,1]$, we have
\[
S_{f}(x)=C_{f}x+O(x^{1/3}(\log x)).
\]
\end{Corollary}

\begin{Remark}
It is believable that one can only get a best possible error term $O(x^{1/3})$ for such type sums. Hence the error term of the square-free case may be improved further by using some deep results related to exponential sums.
%
\end{Remark}

\section{Proof of Theorem \ref{th1}}
Let
\[
x^{\varepsilon}\leq N\leq x^{1-\varepsilon}
\]
be a parameter to be chosen later.
We can write
\[
S_{f}(x):=S_{f,1}(x,N)+S_{f,2}(x,N),
\]
where
\[
S_{f,1}=\sum_{n\leq N}f\left(\left[\frac{x}{n}\right]\right)
\]
and
\[
S_{f,2}=\sum_{N<n\leq x}f\left(\left[\frac{x}{n}\right]\right).
\]
By the assumption in Theorem \ref{th1}, one has
\begin{align}\label{0620}
|f(n)|^{2}\leq \sum_{n\leq x}|f(n)|^{2}\ll x^{\eta}.
\end{align}
Then we have
\[
|f(n)|\ll x^{\eta/2}.
\]
Then  we can obtain that
\begin{align*}
S_{f,1}=\sum_{n\leq N}f\left(\left[\frac{x}{n}\right]\right)
&=\sum_{n\leq N}(x/n)^{\eta/2}
\\& \ll x^{\eta/2}N^{1-\eta/2}\\
&\ll x^{(2+3\eta)/8},
\end{align*}
where $N=x^{1/4}.$

Note that by Cauchy's inequality,  the estimates (\ref{0620}) implies that
\[
\sum_{n\leq x}|f(n)|\ll x^{(1+\eta)/2}.
\]
For $\eta\in(0,2),$ this gives that
\[
\sum_{n\leq N}\frac{f(n)}{n(n+1)}\ll 1.
\]
Hence by using Cauchy's inequality and (\ref{ba0}), we can get
\begin{align*}
S_{f,2}&=\sum_{N<n\leq x}f\left(\left[\frac{x}{n}\right]\right)
\\&
=\sum_{d\leq x/N}f(d)\sum_{x/(d+1)<n\leq x/d}1
\\&
=\sum_{d\leq x/N}f(d)\left(\frac{x}{d}-\frac{x}{d+1}
+\left\{\frac{x}{d+1}\right\}-
\left\{\frac{x}{d}\right\}\right)
\\&
=x\sum_{d=1}^{\infty}\frac{f(d)}{d(d+1)}
-x\sum_{d>x/N}\frac{f(n)}{d(d+1)}
\\&+\sum_{d\leq x/N}|f(d)|
\left(\left\{\frac{x}{d+1}\right\}-
\left\{\frac{x}{d}\right\}\right)
\\&
\ll x\sum_{d=1}^{\infty}\frac{f(d)}{d(d+1)}
+x\sum_{d>x/N}\frac{f(n)}{d(d+1)}
\\&+O\left(\sum_{d\leq x/N}|f(d)|^{2}
\sum_{d\leq x/N}\left(\left\{\frac{x}{d+1}\right\}-
\left\{\frac{x}{d}\right\}\right)^{2}\right)^{1/2}
\\&=x\sum_{d=1}^{\infty}\frac{f(d)}{d(d+1)}
+O\left(\frac{x^{1+\frac{1+\eta}{2}-2}}
{N^{\frac{1+\eta}{2}-2}}
+x^{1/4}\left(\frac{x}{N}\right)^{\frac{\eta}{2}}\right).
\end{align*}
 Observing that $N=x^{1/4},$ we  obtain the desired result.

\section{Proof of Theorem \ref{th2}}
We will start  the proof for Theorem \ref{th2}  with some necessary lemmas.
The following  lemma can be seen in Theorem A.6 in \cite{GK}.
Let $\psi(t)=t-[t]-1/2$ for $t\in\mathbb{R}$ and $\delta\geq0.$ Define
\[
\mathcal{G}(x,D):=\sum_{D<d\leq 2D}f(d)\psi\left(\frac{x}{d+\delta}\right).
\]
\begin{Lemma}\label{z3}
For $0<|t|<1,$ let
$$W(t) = \pi t(1-|t|)\cot\pi t + |t|.$$ For $x\in\mathbb{R},$ $H\geq1,$ we define
$$\psi^{*}(x)=\sum_{1\leq |h|\leq H}(2\pi ih)^{-1}W\left(\frac{h}{H+1}\right)e(hx)$$
and
\[
R_{H}(x)=\frac{1}{2H+2}\sum_{|h|\leq H}\left(1-\frac{|h|}{H+1}\right)e(hx).
\]
Then $\delta(x)$ is non-negative, and we have
$$|\psi^{*}(x)-\psi(x)|\leq
R_{H}(x).$$
\end{Lemma}

In fact, by using some ideas related to Lemma 4.3 in \cite{GK} or Corollary 6.7 in \cite{BB}, we can obtain  the following lemma.
\begin{Lemma}\label{wu}
Let $\delta\geq0$ be a fixed constant. Let $f$ be a positive real-valued arithmetic function such that  $f(n)=\sum_{d|n}g(d)$ and for any sufficiently large $x$
\[
\sum_{n\leq x}|g(n)|\ll x^{\lambda/(1+\kappa)}.
\]
 We have
 \[
\mathcal{G}(x,D)\ll   \left(x^{\kappa}D^{-\kappa+\lambda}
\right)^{1/(1+\kappa)}+x^{-1} D^{2} ,
\]
uniformly for $1\leq D\leq x,$ where $(\kappa,\lambda)$ is an exponent pair.
\end{Lemma}
\begin{Remark}
If we choose
\[
f(n)=\frac{\phi(n)}{n},
\]
in \cite{Wu},
 it is proved that

{\slshape {\bf Theorem A.} (Wu \cite{Wu})
Let $\delta\geq0$ be a fixed constant. We have
 \[
\mathcal{G}(x,D)\ll   \left(x^{\kappa}D^{-\kappa+\lambda}\right)^{1/(1+\kappa)}
+x^{\kappa}D^{-2\kappa+\lambda}(\log x) +x^{-1}D^{2} ,
\]
uniformly for $1\leq D\leq x,$ where $(\kappa,\lambda)$ is an exponent pair. Furthermore, if $(\kappa,\lambda)\neq(1/2,1/2),$ then the factor $\log x$ can be omitted.}

Our result eliminates one   term and the possible $\log x$ in \cite{Wu}.

\end{Remark}

\begin{proof}
Now we will prove Lemma \ref{wu}. The proof relies on the relation
\[
f(d)=\sum_{mn=d}g(m).
\]
Hence one can write
\[
\mathcal{G}(x,D)=\sum_{m\leq 2D}g(m)
\sum_{D/m<n\leq 2D/m}\psi\left(\frac{x}{mn+\delta}\right).
\]
If $D\leq 100x^{\frac{\kappa}{1+2\kappa-\lambda}}
m^{1-\lambda+\kappa}$, then we have
\[
 \left(x^{\kappa}
D^{\lambda-\kappa}\right)^{1/(1+\kappa)}
\gg D \gg\mathcal{G}(x,D).
\]
Hence, we may always assume that $D>100 x^{\frac{\kappa}{1+2\kappa-\lambda}}m^{1-\lambda+\kappa}$.
By Lemma  \ref{z3}, for $x\geq1$ and $H\geq1$, we have
\[
\mathcal{G}(x,D)\ll \sum_{m\leq 2D}|g(m)|\left(
\frac{D/m}{H}+\sum_{h\leq H}\frac1h\left|\sum_{D/m<n\leq 2D/m}e\left(\frac{hx}{mn+\delta}\right)\right|
\right)
\]
with $H\geq 1.$
We also need
the following  well-known lemma (for example, one can refer to page 441 of \cite{BB} or page 34 of  \cite{GK}).
\begin{Lemma}\label{z2}
Let $s^{(k)}(x)\asymp YX^{1-k}$ for $1< X\leq x\leq 2X$ and $k=1,2,\cdots.$ Then one has
\[
\sum_{X<n\leq 2X}e(s(n))\ll Y^{\kappa}X^{\lambda}+Y^{-1},
\]
where $(\kappa,\lambda)\neq(0,1)$ is any exponent pair.
\end{Lemma}

By Lemma \ref{z2}, we have
\begin{align*}
\mathfrak{G}(x,D)&\ll\sum_{m\leq 2D}|g(m)|\left(\frac{D/m}{H}+\sum_{h\leq H}\frac1h\left(\frac{hx}{D^{2}/m}\right)^{\kappa}
\left(D/m\right)^{\lambda}+\frac{D^{2}/m}{hx}\right)\\
&\ll \sum_{m\leq 2D}|g(m)|\left((D/m)H^{-1}+x^{\kappa}
H^{\kappa}(D^{2}/m)^{-\kappa}(D/m)^{\lambda}
+x^{-1}(D^{2}/m)\right).
\end{align*}
Choosing
$$H=\lfloor D^{\frac{1+2\kappa-\lambda}{1+\kappa}}
x^{-\frac{\kappa}{1+\kappa}}/m^{1-\lambda+\kappa}\rfloor$$
 gives the  following result.

\begin{Lemma}[Stucky]\label{St}
Let $\delta\geq0$ be a fixed constant. We have
 \[
\mathcal{G}(x,D)\ll \sum_{d\leq 2D}|g(d)|\left( \left(x^{\kappa}D^{-\kappa+\lambda}d^{-\lambda}
\right)^{1/(1+\kappa)}+x^{-1}d^{-1}D^{2}\right),
\]
uniformly for $1\leq D\leq x,$ where $(\kappa,\lambda)\neq(0,1)$ is an exponent pair.

\begin{Remark}
This lemma can be obtained by  equation (4.2) of \cite{St} by partial summation and the assumption about exponent pairs. As Stucky sketched the proof of equation (4.2) in \cite{St},  here we give a supplement for this.
\end{Remark}
\end{Lemma}
Then the desired conclusion can be obtained by Lemma \ref{St} and the assumption that
\[
\sum_{n\leq x}|g(n)|\ll x^{\lambda/(1+\kappa)}.
\]
\end{proof}

Now we begin the proof of Theorem \ref{th2}. Let
\[
x^{\varepsilon}\leq N\leq x^{1-\varepsilon}
\]
be a parameter to be chosen later.
We can write
\[
S_{f}(x):=S_{f,1} +S_{f,2} ,
\]
where
\[
S_{f,1}=\sum_{n\leq N}f\left(\left[\frac{x}{n}\right]\right)
\]
and
\[
S_{f,2}=\sum_{N<n\leq x}f\left(\left[\frac{x}{n}\right]\right).
\]
Obviously, by the assumption  that $f(n)\ll n^{\varepsilon},$ we have
\begin{align*}
S_{f,1}=\sum_{n\leq N}f\left(\left[\frac{x}{n}\right]\right)
&=\sum_{n\leq N}(x/n)^{\varepsilon}
\\
&\ll N^{1+\varepsilon}
\\ &\ll x^{\lambda/(1+\lambda)+\varepsilon},
\end{align*}
where $N=x^{\lambda/(1+\lambda)}.$

As to $S_{f,2},$ firstly, by the assuming    that  $f(n)=\sum_{d|n}g(d)$ and
\[
\sum_{n\leq x}|g(n)|\ll x^{\lambda/(1+\kappa)},
\]
  we have
\[
\sum_{n\leq x}|f(n)| \ll x .
\]
Hence we can get
\begin{align*}
S_{f,2}&=\sum_{N<n\leq x}f\left(\left[\frac{x}{n}\right]\right)
\\&
=\sum_{d\leq x/N}f(d)\sum_{x/(d+1)<n\leq x/d}1
\\&
=\sum_{d\leq x/N}f(d)\left(\frac{x}{d}-\frac{x}{d+1}
-\psi\left(\frac{x}{d}\right)
+\psi\left(\frac{x}{d+1}\right)\right)
\\&=x\sum_{d=1}^{\infty}\frac{f(d)}{d(d+1)}
+O\left(
N^{1+\varepsilon}
\right)+O\left(\sum_{d\leq x/N}f(d)\psi\left(\frac{x}{d+\delta}
\right)\right).
\end{align*}
Let $x/N=2^{k}$. Then we have
\begin{align*}
\sum_{d\leq x/N}f(d)\psi\left(\frac{x}{d+\delta}
\right)&=\sum_{D=2^{l},\ 0\leq l<k}
\mathcal{G}(x,D)+O(1)\\
&\ll 1+ \sum_{D=2^{l},\ 0\leq l<k}
\left(
x^{\kappa/(1+\kappa)}
D^{(\lambda-\kappa)/(1+\kappa)}
+D^{2}/x\right),
\end{align*}
where we have used Lemma \ref{wu} and the assumption that
 \[
\sum_{n\leq x}|g(n)|\ll x^{\lambda/(1+\kappa)}.
\]

 Then by the estimates of $S_{f,1},$ $S_{f,2}$,
and choosing $N=x^{\lambda/(1+\lambda)},$ we get
\begin{align*}
S_{f}(x)
 =x\sum_{d=1}^{\infty}\frac{f(d)}{d(d+1)}
+O\left(
 x^{\lambda/(1+\lambda)+\varepsilon}
 +(\log x)^{\alpha}x^{1-2\lambda/(1+\lambda)}\right).
\end{align*}
Recall the fact such that $\lambda\geq 1/2,$ we can finish the proof of Theorem \ref{th2} for $|f(n)|\ll n^{\varepsilon}$.
For $|f(n)|\ll 1,$ similar arguments give that
\begin{align*}
S_{f}(x)
 =x\sum_{d=1}^{\infty}\frac{f(d)}{d(d+1)}
+O\left(
 x^{\lambda/(1+\lambda)}
 +(\log x)^{\alpha}x^{1-2\lambda/(1+\lambda)}\right).
\end{align*}

\section{Remarks on Corollary \ref{co}}
Let $\mu_{k}(n)$ be the indicator function of the $k$-free numbers.
Then we have
\[
\mu_{k}(n)=\sum_{d|n}g(d)
\]
with
 \begin{align*}
g(d)=
 \begin{cases}
\mu(l)\ \ &\textup{if}\ d=l^{k},\\
\\
0\ \ &\textup{otherwise}.
\end{cases}
 \end{align*}
 Then we have
 \begin{align*}
\sum_{d\leq x}g(d)=
 \begin{cases}
 O\left(x^{1/2}\right)\ \ &\textup{if}\ k=2,\\
\\
 O\left(x^{1/3 }\right)\ \ &\textup{if}\ k\geq3.
\end{cases}
 \end{align*}
Choosing $(\kappa,\lambda)=(1/2,1/2)$ in Theorem \ref{th2} gives the case of $k\geq 3$ in Corollary \ref{co}.
In fact, as $\mu_{3}(n)\ll 1,$ we have
\begin{align*}
 S_{\mu_{3}}(x)=
\sum_{n=1}^{\infty}
\frac{\mu_{3}(n)}{n(n+1)}x+O\left(x^{1/3}(\log x)^{2}\right).
 \end{align*}
 And one $\log$ can be cancelled  by partial summation for $k>3.$ Hence for $k>3,$ we have
 \begin{align*}
 S_{\mu_{k}}(x)=
\sum_{n=1}^{\infty}
\frac{\mu_{k}(n)}{n(n+1)}x+O\left(x^{1/3}(\log x)\right).
 \end{align*}
Choosing $(\kappa,\lambda)=BABAAB(0,1)=(4/18,11/18)$ in Theorem \ref{th2}, we can obtain Corollary \ref{co} for $k=2.$ However, for $k=2,$ the $\log$ term in the error term is   not really important since the exponent  is probably not the best possible one.

In fact, inspired by the recent work of the author \cite{zw111}, we can obtain a slightly better result for the case of $k=2$ in Corollary \ref{co}.
Assuming that $(a,1/2+a)$  is an exponent  pair, then  by
$$BABA(a,1/2+a)=\left(\frac{2a+1}{6a+5},
\frac{4a+3}{6a+5}\right),$$
we can get
  \begin{align*}
 S_{\mu_{2}}(x)=
\sum_{n=1}^{\infty}
\frac{\mu_{2}(n)}{n(n+1)}x
+O\left(x^{\vartheta(a)+\varepsilon}\right)
 \end{align*}
 such that
 \[
 \vartheta(a)=\frac{4a+3}{10a+8}.
 \]
Choosing $a=1/6,$   we have
  \begin{align*}
 S_{\mu_{2}}(x)=
\sum_{n=1}^{\infty}
\frac{\mu_{2}(n)}{n(n+1)}x
+O\left(x^{11/29+\varepsilon}\right).
 \end{align*}
This gives the result of $k=2$.
By the work of Bourgain \cite{Bo}, we can choose $a=13/84+\varepsilon.$ Then we have
  \begin{align*}
 S_{\mu_{2}}(x)=
\sum_{n=1}^{\infty}
\frac{\mu_{2}(n)}{n(n+1)}x
+O\left(x^{152/401+\varepsilon}\right),
 \end{align*}
 which gives an improvement of $11/29.$

\bigskip
{\bf Acknowledgements}
I am deeply grateful to the referee(s) for carefully reading the manuscript and making useful suggestions and significant corrections.

\address{Wei Zhang\\ School of Mathematics and Statistics\\
               Henan University\\
               Kaifeng  475004, Henan\\
               China}
\email{zhangweimath@126.com}

\end{document}